\documentclass[10pt]{article}
\usepackage[english]{babel}
\usepackage{latexsym}
\usepackage{color}
\usepackage{amsmath}
\usepackage[latin1]{inputenc}
\usepackage{amsfonts}
\usepackage{amssymb}
\usepackage{mathrsfs}
\usepackage{eucal}
\usepackage{exscale}
\usepackage{makeidx}
\usepackage{makeidx}
\usepackage{graphicx}
\usepackage[T1]{fontenc}
\usepackage{amsmath,amssymb}
\usepackage[all]{xy}
\newcommand{\qed}{{$\square$}\medbreak}
\newtheorem{thm}{Theorem}%[subsection]
\newtheorem{rem}{Remark}%[subsection]
\newtheorem{df}{Definition}

\newtheorem{cor}{Corollary}
\newtheorem{lem}{Lemma}

\newtheorem{exemple}{Example}

\newcommand{\preuve}{\indent {\it Proof.}\hspace{4mm}}

\title{Positive Strong Amalgamations}
\author{Mohammed Belkasmi}

\begin{document}
\maketitle

\begin{abstract} 
In  this paper, we present the notions 
of positively complete theory and general forms 
of amalgamation in the framework of positive logic.
We explore the fundamental properties of 
positively complete theories and study the behaviour 
of companion theories by a change of constants in the language. 
Moreover, we present  a general form of amalgamation and 
discuss some forms of strong amalgamation.
 \end{abstract} 
\section{Positively complete theories}
\subsection{Positive logic}
The positive logic  in its present form
was introduced
 by Ben Yacccov and Poizat  \cite{begnacpoizat},
  following the line of research of   
  Hrushovski \cite{rob} and Pillay \cite{pillay}.
 
 Broadly speaking, the positive model theory is
 considered as a part of 
 the eastern model theory introduced by 
 Abraham Robinson, wich is  concerned 
 essentially with  the study of
existentially closed models and  model-complete theories 
 in the context of  incomplete inductive theories. 
 The main tools in the study of incomplete inductive theories are:
  embedding, existential formulas and inductive sentences. 
   Keep in consideration homomorphisms 
 and positive formulas, the positive logic offers 
 a wider and simpler framework as companed to
  the eastern model theory.

In this subsection we summarise  the basic concepts 
of positive logic which will be used 
throughout the paper.

Let $L$ be a first order language. we stipulate
that $L$ includes the symbol of
equality and the constant $\perp$ denoting the antilogy.\\
The quantifier-free positive formulas are built from atomics 
by using the connectives $\wedge$ and
$\vee$. The positive formulas are of the form: 
$\exists\bar x\varphi(\bar x, \bar y)$,
 where $\varphi$  is quantifier-free positive formula.\\
A sentence is said to be h-inductive,  if it is a finite
conjunction of sentences of the form:
$$\forall\bar x\\ (\varphi(\bar x)\rightarrow 
\psi(\bar x))$$
where $\varphi$ and $\psi$ are  positive formulas.\\
The h-universal sentences are the sentences that can be written as negation of a positive sentence.

Let $A$ and $B$ be two $L$-structures and 
 $f$ a mapping from $A$ into $B$. $f$ is said to be
 \begin{itemize}
 \item  a homomorphism, if for every  tuple $\bar a$ from  $A$ ($\bar a\in A$ by abuse of notation)
and for every atomic formula $\phi$, $A\models\phi(\bar a)$
implies $B\models\phi(f(\bar a))$. In this  case we say 
that  $B$ is a continuation of $A$.
\item an embedding, if $f$ is a  homomorphism 
 such that  for every
atomic formula $\phi$; $A\models\phi(\bar a)$
if and only if $B\models\phi(f(\bar a))$.
\item  an immersion whenever  $\bar a\in A$ and 
 $f(\bar a)$  satisfy the same $L$-positive formulas, for every 
 $\bar a\in A$.
 \end{itemize}
 Given $A$ be a $L$-structure, we let 
 $L(A)$ be the language obtained from $L$ by adjoining 
the element of $A$ as constants. 
We denote by $Diag(A)$
(resp. $Diag^+(A)$) the 
set of   atomic and negated atomic 
(resp. positive quantifier-free) sentences  satisfied by
 $A$ over the 
language $L(A)$.\\
We denote by $Diag^{+\star}(A)$ the 
set of    positive sentences  true in $A$ over the 
language $L$.
\begin{df}
A $L$-structure (resp. model) $M$ of a class $\Gamma$   (resp. 
of an h-inductive $L$-theory $T$)
is said to be pc  in $\Gamma$ (of $T$)
 if every homomorphism from $M$ into 
a member
of $\Gamma$ (resp. into a model of $T$)
is an immersion. 
\end{df}
In \cite{begnacpoizat} it is shown that the  pc models exist for 
any consistent h-inductive theory. We have the 
useful following fact.
\begin{lem}[{\cite[Th\'eor\`eme 1, lemme 12] {begnacpoizat}}]\label{pecconti}:
\begin{itemize}
\item Every member of an $h$-inductive class is continued 
in a pc member of the class.
\item The class of pc models of an $h$-inductive theory $T$ is
 $h$-inductive.
\end{itemize}
\end{lem}
\subsection{Positively complete and T-complete theories}
 \begin{df}
Two h-inductive theories are said to be companion 
 if they 
have the same pc models.
\end{df}
 Every $h$-inductive 
theory $T$ has a maximal companion denoted $T_k(T)$,  called
the Kaiser's hull of $T$.  $T_k(T)$
is the $h$-inductive theory of the pc models of $T$.
 Likewise,
$T$ has a minimal companion denoted $T_u(T)$, formed by its
$h$-universal consequences sentences. 

Note that if $T'$ is an h-inductive theory such that 
$T_u(T)\subseteq T'\subseteq T_k(T)$ then $T'$ and 
$T$ are companion theories.
\begin{df}
Let $T$ be an h-inductive theory.
\begin{itemize}
 \item $T$ is said to be model-complete 
if  every model of $T$ is a pc model of $T$.
\item We say that $T$ has a model-companion whenever $T_k(T)$ is 
model-complete.
\end{itemize}
 \end{df}
 Let $A$ be a $L$-structure and 
 $B$ a subset of $A$. 
 We shall use the following notations:
\begin{itemize}
 \item $T_i(A)$ (resp. $T_u(A)$) denote the set of h-inductive
 (resp. h-universal ) 
 $L(A)$-sentences satisfied by $A$.
\item $T_i^\star(A)$ (resp. $T_u^\star(A)$) 
denote the set of h-inductive (resp. h-universal) $L$-sentences 
 satisfied by $A$.
 \item $T_k(A)$ (resp. $T_k^\star(A)$) denote the Kaiser's
hull of $T_i(A)$
(resp. of  $T_i^\star(A)$).
\item $T_i(A| B)$ (resp. $T_u(A| B)$) denote the set of h-inductive
 (resp. h-universal ) 
 $L(B)$-sentences satisfied by $A$.
\end{itemize}
\begin{df}
 Let $A$ and $B$ two $L$-structures and $f$ a homomorphism 
 from $A$ into $B$. $f$ is said to be a strong immersion 
 if $B$ is a model of $T_i(A)$ in the language $L(A)$.
 \end{df}
\begin{df}
\begin{itemize}
\item An h-inductive theory $T$ is said to be positively complete
(or it has 
the joint continuation (in short JC) property) if  
 any two models of $T$ have a common continuation.
 \item Let $T_1, T_2$ and $T$ three h-inductive
 $L$-theories.  $T_1$ and $T_2$ are said 
 to be $T$-complete if for every models $A$ of $T_1$
 and $B$ of $T_2$, there is $C$ a common continuation 
 of $A$ and $B$ such that $C\vdash T$.
\end{itemize}
 
\end{df}
The following remark lists some simple 
 properties which will be 
useful in the rest of the paper.
\begin{rem}\label{remkpropTiTv}
Let $A$ and $B$ two $L$-structures and $T$ an h-inductive $L$-theory.
\begin{enumerate}
\item  $ T_u(A)$ and $ Diag^+(A)$ 
(resp $ T_u^\star(A)$ and $Diag^{+\star}(A)$)
 are subsets of $T_i(A)$ (resp. $T_i^\star(A)$).
\item  $A$ is a pc model of $T_i(A)$, and 
$T_i(A)= T_k(A)$.
\item $T_u(T)$ (resp. $T_u(A)$) is the h-universal part of 
$T_k(T)$ (resp. $T_i(A)$). The same is true for 
  $T_u^\star(A)$ and $T_i^\star(A)$.
\item $T_i(A)\subseteq T_i(B)\Rightarrow 
T_u(A)\subseteq T_u(B)$.
\item $T_i^\star(A)\subseteq T_i^\star(B)\Rightarrow 
T_u^\star(A)\subseteq T_u^\star(B)$.
\item If $T$ is positively complete and $A$ a pc model of $T$, then 
$T_k(T)= T_i^\star(A)$ and $T_u(T)=T_u^\star(A)$.
\item If $A$ and $B$ are pc models of $T$ and $B$ is a continuation of 
$A$ then $T_i^\star(A)= T_i^\star(B)$.
 
\item If $A$ is continued in $B$ then $T_u^\star(B)\subseteq T_u^\star(A)$.
\item If $A$ is immersed in $B$ then $T_u^\star(A)= T_u^\star(B)$ and 
$T_i^\star(B)\subseteq T_i^\star(A)$.
\item $T_u^\star(A)=\{\neg\exists\bar x\varphi(\bar x)\ \mid
\exists\bar x\varphi(\bar x)\notin Diag^{+\star}(A)\}$.
\item $Diag^{+\star}(A)\subseteq Diag^{+\star}(B)\Leftrightarrow
T_u^\star(B)\subseteq T_u^\star(A)$.
\item If 
$T_u^\star(A)\subseteq T_u^\star(B)$ 
(resp. $T_i^\star(A)\subseteq T_i^\star(B)$), then 
$T\cup Diag^+(A)\cup Diag^+(B)$  is consistent in
the language $L(A\cup B)$. 
\item If $T_u^\star(A)\subseteq T_u^\star(B)$ then 
$Diag^+(A)\cup Diag^+(B)$  is consistent over 
the language $L(A\cup B)$.
\item For every pc models $A$ and $B$ of $T$, 
if $T_u^\star(A)= T_u^\star(B)$ then 
$T_i^\star(A)= T_i^\star(B)$.
\item $T_1$ and  $T_2$ are $T$-complete if and only if 
for every $A\vdash T_1$ and $B\vdash T_2$,  
$Diag^+(A)\cup Diag^+(B)\cup T$ is 
$L(A\cup B)$-consistent.
\item Let $(T_1, T_2)$ be a pair of $T$-complete theories.
For every $T_1', T_2'$ and $T'$ companion theories
of $T_1, T_2$ and $T$ respectively, the 
pair  $(T_1', T_2')$ is $T'$-complete.
\end{enumerate}
\end{rem}

\begin{lem}\label{latticeoftheory}
 Let  $A$ be a pc model of an h-inductive $L$-theory $T$, then
 \begin{enumerate}
 \item $T_u^\star(A)$  is minimal 
 in the set 
$\{T_u^\star(B)\mid B\models T\}$. 
\item $T_i^\star(A)$ is 
 maximal in the set 
$\{T_i^\star(B)\mid B\models T\}$. 
\end{enumerate}  
\end{lem}
\preuve 
\begin{enumerate}
 \item  Let $B$ a model of $T$ such that 
$T_u^\star(B)\subseteq T_u^\star(A)$. By 
the property $13$ of the Remark \ref{remkpropTiTv}, there exists 
$C$ a model of $T$ that is a common continuation of 
$A$ and $B$. Given that $A$ is a pc model, from the properties 
$8$ and $9$ of the Remark \ref{remkpropTiTv}
 it follows that 
$$T_u^\star(A)=T_u^\star(C)\subseteq T_u^\star(B).$$ 

\item  Let $B$ a model of $T$ such that 
$T_i^\star(A)\subseteq T_i^\star(B)$. We claim that 
$Diag^+(A)\cup T_i(B)$ is consistent in the language 
$L(A\cup B)$. Indeed, if not, by compactness there exists
$\psi(\bar a)\in Diag^+(A)$ such that 
$T_i(B)\vDash\neg\exists\bar x \psi(\bar x)$. Given that 
$T_i^\star(B)$ is the part of $T_i(B)$ without parameters of 
$B$, then $T_i^\star(B)\vDash\neg\exists\bar x \psi(\bar x)$.
On the other hand since  
$$\exists\bar x \psi(\bar x)\in Diag^{+\star}(A)\subset
 T_i^\star(A)\subseteq T_i^\star(B),$$
 we obtain a contradiction. Thereby $Diag^+(A)\cup T_i(B)$ is consistent in the language 
$L(A\cup B)$, which implies the existence of a model $D$ of 
$T_i(B)$ in the language $L(A\cup B)$, such that
 \[
\xymatrix{
    A \ar[r]^{f} & {D} & {B} \ar[l]_{g}.     
  }
\]
where $f$ is an homomorphism and $g$ an immersion.\\
Given that $D$ is also a model of $T$ and $A$ pc model of $T$, then 
$f$ is an immersion. By the property $9$ of the 
remark \ref{remkpropTiTv} we obtain 
$$T_i^\star (B)\subseteq T_i^\star(D)\subseteq T_i^\star(A)
\subseteq T_i^\star(B).$$
\end{enumerate}
\begin{lem} Let 
$T_1, T_2$ and $T$ three h-inductive $L$-theories.
$T_1$ and $T_2$ are $T$-complete if and only if 
one of the following holds:
\begin{enumerate}
\item For every free-quantifier positive formulas
 $\varphi(\bar x)$ and $\psi(\bar y)$,
If $T\vdash\neg\exists\bar x\varphi(\bar x)
\vee\neg\exists\bar y\psi(\bar y)$ then 
$T_1\vdash\neg\exists\bar x\varphi(\bar x)$ or 
$T_2\vdash\neg\exists\bar y\psi(\bar y)$.
\item $T_u(T)\subseteq T_u(T_1)\cap T_u(T_2)$.
\end{enumerate}
\end{lem}
\preuve 
\begin{enumerate}
\item Suppose that  $T_1, T_2$ and $T$ satisfy the 
following:
\begin{itemize}
\item $T_1$ and $T_2$ are $T$-complete.
\item $T\vdash\neg\exists\bar x\varphi(\bar x)
\vee\neg\exists\bar y\psi(\bar y)$.
\item $T_1\nvdash\neg\exists\bar x\varphi(\bar x)$.
\item $T_2\nvdash\neg\exists\bar y\psi(\bar y)$.
\end{itemize}
 Then 
there are $A$ and $B$ models of $T_1$ and 
$T_2$ respectively, such that 
$A\models \varphi(\bar a)$ for some $\bar a\in A$,
and $B\models\psi(\bar b)$ for some $\bar b\in B$.
Let $C$ be a common continuation of $A$ and $B$ that 
is a model of $T$, then 
$C\models\exists\bar x\varphi(\bar x)\wedge
\exists\bar y\psi(\bar y)$, contradiction.\\
Conversely, suppose that $T_1, T_2$ and $T$ satisfy 
the property $1$. Let $A$ and $B$ models of $T_1$ and
$T_2$ respectively. We claim that 
$Diag^+(A)\cup Diag^+(B)\cup T$ is 
$L(A\cup B)$-consistent. If not, there are 
$\varphi(\bar a)\in Diag^+(A)$ and 
$\psi(\bar b)\in Diag^+(B)$ such that 
$T\vdash\neg(\exists\bar x\varphi(\bar x)\wedge
\exists\bar y\psi(\bar y))$. Thereby 
$T_1\nvdash\neg\exists\bar x\varphi(\bar x)$
and $T_2\nvdash\neg\exists\bar y\psi(\bar y)$, 
contradiction.
\item Suppose that $T_1$ and $T_2$ are $T$-complete,
it is clear that 
$T_u(T)\subseteq T_u(T_1)\cap T_u(T_2)$.\\
The proof of the other direction  is the same 
as the second part of the proof of $1$ applies 
at the theories 
$T_u(T_1),  T_u(T_2)$ and $ T_u(T)$, knowing that
$T_1$ and $T_2$ are $T$-complete if and only if
$T_u(T_1)$ and  $T_u(T_2)$ are $ T_u(T)$-complete.

\end{enumerate}

\begin{lem}\label{maxcompleteth}
An h-inductive $T$ theory has the JC property if and only if one of the following holds:
\begin{enumerate}

\item For any free-quantifier positive formulas $\varphi(\bar x)$ and $\psi(\bar y)$,
if $T\vdash\neg\exists\bar x\varphi(\bar x)
\vee\neg\exists\bar y\psi(\bar y)$ then 
$T\vdash\neg\exists\bar x\varphi(\bar x)$ or 
$T\vdash\neg\exists\bar y\psi(\bar y)$.
\item For every positive formulas 
$\varphi(\bar x)$ and $\psi(\bar y)$, if 
$T\cup\{\varphi(\bar x)\}$ and 
$T\cup\{\psi(\bar y)\}$ are consistent sets then 
$T\cup\{\varphi(\bar x), \psi(\bar y)\}$ is a 
consistent set.
\item  $T_u(T)= T_u^\star(A)$ for some model $A$ of $T$. 
\item  $T_k(T)= T_i^\star(A)$ for some model $A$ of $T$. 
\item For every pc models $A$ and $B$ of $T$ we have 
$T_u^\star(A)= T_u^\star(B)$.
%\item If $T\subseteq T'$ and $T$ is positively complete, then so is 
%$T'$. 
\end{enumerate}
\end{lem}
\preuve 
\begin{enumerate}
\item Well known for complete theories in first 
order logic.
\item Clear
\item Let $T$ be an h-inductive theory and $A$ a model 
of $T$ such that $T_u(T)= T_u^\star(A)$.
Let $B$ and $C$  two pc models of $T$. Given that  
$T_u(T)= T_u^\star(A)\subseteq T_u^\star(B)\cap T_u^\star(C)$, by
the minimality of the h-universal theory of the pc models, we obtain 
 $$T_u^\star(A)=T_u^\star(B)=T_u^\star(C).$$
 From the property $13$ of the Remark \ref{remkpropTiTv}, we deduce 
the existence of a common continuation of $B$ and $C$ by a 
model of $T$.
Thereby $T$ is positively complete.

The other direction follows from $6$ of the remark \ref{remkpropTiTv}.
\item Let $A$ be a model of $T$ such that $T_k(T)= T_i^\star(A)$.
Let $B$ and $C$ be two pc models of $T$. Since 
$$ T_i^\star(A)= T_k(T)\subseteq T_i^\star(B)\cap T_i^\star(C)$$
then $ T_u^\star(A)\subseteq T_u^\star(B)\cap T_u^\star(C)$. Given that 
$T_u^\star(B)$ and $T_u^\star(C)$ are minimal, we obtain
$$ T_u^\star(A)= T_u^\star(B)= T_u^\star(C).$$
By   
the property $13$ of the 
Remark \ref{remkpropTiTv}, we get  
 a common continuation of $B$ and $C$ by a 
model of $T$.
Thereby $T$ is positively complete.

The other direction results from the 
property $6$ of the remark \ref{remkpropTiTv}.
\end{enumerate}
\begin{lem}\label{TvTicompanion}
Let $A$ be a $L$-structure. The theories $T_u^\star(A)$ and 
$T_i^\star(A)$ are companion and positively 
completes. 
\end{lem}
\preuve 
 It is clear that every 
 model of $T_i^\star(A)$ is a model of $T_u^\star(A)$.\\ 
 In the next step  we will show that every model of 
 $T_u^\star(A)$
 is continued into a model of $T_i^\star(A)$.  Let 
 $B$ be a model of $T_u^\star(A)$, we claim that  
 $Diag^+(B)\cup T_i^\star(A)$ is  consistent  
 in the language $L(B)$. Indeed, if not, 
  by compactness 
 there exists $\psi( \bar b)\in Diag^+(B)$ such that
 $T_i^\star(A)\vDash \neg\exists\bar x\psi(\bar x)$, then 
 $\neg\exists\bar x\psi(\bar x)\in T_u^\star(A)$.
 Given that $ T_u^\star(A)\subseteq T_u^\star(B)$
 and $\exists\bar x\psi(\bar x)\in Diag^{+\star}(B)$,
 we obtain a contradiction. Thereby 
 $Diag^+(B)\cup T_i^\star(A)$ is  consistent,  so 
 $B$ is continued in a model of $ T_i^\star(A)$.
 
 The second part of the lemma results from 
 the properties $2$ and $3$ of the lemma \ref{maxcompleteth}.\qed
%$T_u^\star(A)= T_u(T_u^\star(A))$. It follows from 
% the property 
%  that $T_u^\star(A)$ is positively complete.\qed
\begin{rem}\label{rem2}
\begin{itemize}
\item We have the same results of the lemma \ref{TvTicompanion} for the theories $T_u(A| B)$ and $T_i(A| B)$, where 
 $B$ is a subset of $A$.
 \item 
 Let $A_e$ be a pc model of an h-inductive theory $T$.
 Let $A$ be a subset of $A_e$. Every pc model of 
 $T_u(A_e| A)$ in the language $L(A)$ is a 
 pc model of $T$ in the language $L$.\\
 Indeed, Let $B_e$ be a pc model of $T_u(A_e| A)$,
 since $T_u(A_e| A)$ is positively complete, 
 there is a common continuation $C$ of 
 $A_e$ and $B_e$ in the language $L(A)$ which in turn 
 can be continued in a pc model $C_e$ of $T$. 
 As $A_e$ is immersed in $C_e$, so $C_e$ is a 
 model of $T_u(A_e| A)$, then $B_e$ is immersed 
 in $C_e$, which implies that $B_e$ is a pc model of 
 $T$.
 \end{itemize}
 \end{rem}
  \begin{lem}\label{lemma5}
  Let $T$ be a positively complete h-inductive 
  $L$-theory and 
  $A_e$ a pc model of $T$ that is also a pc model of 
  an h-inductive $L$-theory $T'$. Then every pc model of $T$
  is a pc model of $T'$, and every pc model 
  of $T'$ that is  
  a model of $T$ is a pc model of $T$.
  \end{lem}
 \preuve 
 Given that $A_e$ is a pc model of $T'$, then 
 $T'\subset T_k(T)= T_i^\star(A_e)$. Let $B$ be a pc model of 
 $T$, since $T_i^\star(A_e)= T_i^\star(B)$ then $B$ is a model 
 of $T'$.\\
  Let $f$ be a homomorphism from $B$ into $B'$ a pc model 
 of $T'$, from  the property $8$ of the remark 
 \ref{remkpropTiTv} we have
 $$T_u^\star(B')\subseteq T_u^\star(B)=T_k(T)= T_u^\star(A_e).$$
 Now by the property $12$ of the remark
  \ref{remkpropTiTv}, we obtain the consistency  of 
 $T'\cup Diag^+(A_e)\cup Diag^+(B')$, which gives 
 the following diagram:
 \[
\xymatrix{
   & A \ar[rd]^{i_m}   \\
     {B} \ar[r]^{f}& B' \ar[r]_{f'} & {C}
  }
\]
where $C$ is a model of $T'$ that we can take it  pc model. We deduce the following equalities:
 $$T_k(T)=T_i^\star(B)= T_i^\star(A_e)= 
 T_i^\star(C)=T_i^\star(B').$$
 Thereby $f$ is an immersion, and $B$ is a pc model of $T'$.
 
 For the second part of the lemma. Let $B_e$ be a pc model of 
 $T'$ such that $B_e\vdash T$, let $f$ be a homomorphism from 
 $B_e$ into a pc model $B$ of $T$. Given that $B$ is also 
 a pc model of $T'$, then $f$ is an immersion, and so 
 $B_e$ is a pc model of $T$.\qed
 \begin{cor}\label{corpcofTi(A)}
 Let $T$ be an h-inductive theory and $A$ a pc model of $T$. 
 Every pc model of the $L$-theory
  $T_i^\star(A)$ is a pc model of $T$, and every 
 pc model of $T$ which is a model of $T_i^\star(A)$ is a pc
 model of $T_i^\star(A)$.
 \end{cor}
 \preuve 
 The corollary follows directly from the fact that 
 $T_i^\star(A)$ is positively complete and $A$ is a common pc 
 of $T$ and $T_i^\star(A)$.\qed 
 \begin{rem}
 \begin{itemize}
 \item  If the 
 language of theory $T'$ in the lemma \ref{lemma5}
 contains the language of $T$, 
 we obtain a similar result.
 We have the  possibility of interpreting a
  pc model of $T$ in the language of $T'$. 

 \item Let $A$ be a subset of $A_e$ a pc 
 model of $T$, let
 $<A>$ be the subset of $A_e$ generated by  $A$. 
 It follows from  the second property 
 of the  remark  \ref{rem2} that every pc 
model of  $T_u(A_e| <A>)$ in the language 
 $L(<A>)$ is a pc model of  $T_u(A_e| A)$ in 
 the language $L(A)$.\\
 Conversely, every pc model of $T_u(A_e| A)$
 is a pc model of $T_u(A_e| <A>)$. Indeed, let 
 $B_e$ be a pc model of $T_u(A_e| A)$, 
 since $A_e$ is a pc model of the positively 
 complete theory $T_u(A_e| A)$, there exist 
a pc model  $C_e$ which is a common continuation
 of $A_e$ and $B_e$. 
 thereby the substructures generated by 
 $A$ in $A_e, B_e$ and $C_e$ are isomorphic, which 
 implies that $B_e$ is a model of $T_u(A_e| <A>)$.
 Therefore $B_e$ is a pc model of $T_u(A_e| <A>)$.
 \end{itemize}

 \end{rem}

\begin{lem}
Let $T$ be a positively complete h-universal theory 
over a language containing a non empty set of 
constants $C$. If for 
every model $A$ of $T$ we have $T_u^\star(A)= T$,
then there exist a model $M$ of $T$ such that the class
of pc models of $T$ is equal to the class of pc models
of $T_u(M)$ in the language $L(M)$.
\end{lem}
\preuve 
Let $A$ be a pc model of $T$, denote by 
$C_A$ the set of interpretation of $C$ in $A$ and by 
$<C_A>$ the substructure of $A$ generated by $C_A$.\\
Given that $T$ is positively complete, then for every 
pc models $A$ and $B$ there exists a pc model 
$C$ of $T$ that is a common continuation of $A$ and 
$B$, thereby the structures $<C_A>, <C_B>$ and 
$<C_C>$ are isomorphic. Consequently, every pc model
$A$ of  $T$ is a pc model of $T_u(<C_A>)$ in the 
language $L(<C_A>)$.\\
Conversely, let $A$ be a pc model of $T$, we will show 
that $<C_A>$ is a pc model of $T$. 
Considering that the elements of  $<C_A>$ 
are the terms of the language $L$ (modulo $T$) 
and $<C_A>$ is a model of $T$, we can 
suppose that $T_u(T)=T_u(<C_A>)$ in the language 
$L(<C_A>)$. Since $<C_A>$ is embedded in $A$
and $A$ is a model of $T_u(<C_A>)$ then 
$<C_A>$ is immersed in $A$, thereby $<C_A>$
is a pc model of $T$.\qed 
The following example list some anomaly situations 
in the positive logic that we will try to deal
by some changes focused on  the language and the theories. 
\begin{exemple}
\begin{enumerate}
\item Let $T_{pos}$ the h-inductive theory of posets 
in the relational language $L=\{\leq\}$.
$T_{pos}$ is positively complete 
and has only one pc model which is 
the trivial structure $(\{x\}, \leq)$.
\item Let $L$ be the language formed by the 
functional symbol $L$.
\begin{enumerate}
\item For every integer $n$, let $T_n$ be the h-inductive theory $\{\exists x\ f^n(x)=x\}$.
For every $n$, the theory $T_n$ is positively complete 
and has only one pc model which is the structure 
$(\{x\}, f)$ such that $f(x)=x$.
\item For every integer $n$, let let $T_n$ be the h-inductive theory $\{\neg\exists x\ f^n(x)=x\}$. 
For every integer $n$,  we can
view the models of $T_n$ as directed graphs 
such that the vertexes of the graph are 
the element of the model, and two 
vertexes $a$ and $b$ are jointed by an 
edge pointed from $a$ into $b$ if $f(a)=b$. 
The theory $T_n$ is 
positively complete and has only one pc formed 
by  the graph that contains for every prime 
$p$ that not divide $n$, 
one  cycles of length $p$ 
\end{enumerate}
\item Let $T_g$ the h-inductive theory of groups 
 in the useul  language $L_g$ of groups. The trivial group
 is the unique pc model of $T_g$.
\end{enumerate}
\end{exemple}

In order to rectify the anomaly observed in the last 
example, we  propose two distinct methods which we will apply to the theory $T_g$.\\
the first method consists  to define positively the inequality by adding a binary relation symbol $R$
to $L_g$ interpreting  by 
$R(a, b)\leftrightarrow a\neq b$. The pc model 
of the new theory so defined are the existentially 
closed groups in the context of logic with negation.\\
The second method consists to discard the trivial 
groups by adding a symbol of constant $a$ to 
$L_g$ and consider the theory $T_g^+= T_g\cup 
\{a\neq e\}$.
In the following we will draw some features of 
the class of pc models of $T_g^+$. Let $G$
be a pc model of $T_g^+$ and $a_G$ the interpretation 
of the new constant in $G$, let $L^+$ the language 
of $T_g^+$.
\begin{itemize}
\item Unlike $T_G$, the theory $T_g^+$ is not 
positively complete. 
\item For every integer $n$ there exists an element 
of $G$ of order $n$ (just embed $G$ into a group 
that satisfy $\exists x\ x^n=e$).
\item The pc models of  $T_g^+$ are simple. Indeed,
suppose that $G$ is not simple,
let $N$ be a normal subgroup of a pc model $G$ 
and $f_N$ the canonical homomorphism  defined 
from $G$ into $G/N$. To make $f_N$ an immersion, it is necessary that $a_G\in N$, thereby $a_G$ must belong 
to every normal subgroup of $G$.\\
Given that every  $L^+$-homomorphism is an homomorphism
of groups, it follows that $f(a_G)=e$ for every 
homomorphism of groups, contradiction.
\item $(G, b_G)$ is a pc model for every $b\neq e$ in 
$G$.
Indeed, if $f$ is an $L_g(b_G)$-homomorphism  then $f$ 
is an $L_g(a_G)$-homomorphism (since $f(a_G)\neq e$),
thereby $f$ is an $L_g(b_G)$-immersion.
\end{itemize}

\section{General forms of amalgamation }
For our needs, we adopt the following notations:\\
Let $A$ and $B$ are two $L$-structures and $f$ a mapping 
from $A$ into $B$. We say that $f$ is:
\begin{itemize}
\item $h-hom$ if $f$ is a homomorphism.
\item $e-hom$ if $f$ is an embedding.
\item $i-hom$ if $f$ is an immersion.
\item $s-hom$ if $f$ is a strong immersion.
\end{itemize}

 \begin{df}\label{dfamalgamation}
 Let $\Gamma$ be a class of $L$-structures and $A$ a member of 
 $\Gamma$. Let $\alpha, \beta, \gamma, \delta$  members of  
 $\{h-hom, e-hom, i-hom, s-hom\}$.
  We say that $A$ is an 
  $[\alpha, \beta, \gamma, \delta]$-amalgamation basis of 
 $\Gamma$
if  for all $B, C$ in $\Gamma$,  $f$ a 
$\alpha-hom$ 
 from $A$ into $B$ and $g$ a $\beta-hom$ from 
$A$ into $C$,
there exist 
$D\in\Gamma$, $f'$ a $\gamma-hom$ from $B$ into $D$
and  $g'$ a $\delta-hom$ from $C$ into $D$
 such that
the following diagram commutes:
\[
\xymatrix{
    A \ar[r]^{f} \ar[d]_{g} & {B} \ar[d]^{g'} \\
    C \ar[r]_{f'} & {D}
  }
\]
 We say that $\Gamma$ has the
   $[\alpha, \beta, \gamma, \delta]$-amalgamation 
property if every element of $\Gamma$ is an 
$[\alpha, \beta, \gamma, \delta]$-amalgamation basis of 
$\Gamma$.\\
 We say that $A$ is an 
$[\alpha, \beta]$-asymmetric amalgamation basis of $\Gamma$,
if  $A$ is $[\alpha, \beta, \alpha, \beta]$-amalgamation basis of 
$\Gamma$.\\
 We say that $A$ is  
$[\alpha, \gamma]$-pregeneric amalgamation basis of $\Gamma$,
if $A$ is $[\alpha, \alpha, \beta, \beta]$-amalgamation basis of 
$\Gamma$.\\
 We say that $A$ is an 
$[\alpha]$- amalgamation basis of $\Gamma$,
if$A$ is $[\alpha, \alpha, \alpha, \alpha]$-amalgamation basis of 
$\Gamma$.\\
We say that $A$ is $[\alpha, \gamma]$-strong amalgamation basis 
of $\Gamma$, if $A$ is  $[\alpha, \gamma]$-pregeneric
 amalgamation basis of $\Gamma$, 
 and for every $B, C$ members of $\Gamma$; if 
$A$ is continued into $B$ and $C$ by $\alpha-hom$, then there  exist
$D\in \Gamma$,  $f'$ and $g'$ $\gamma-hom$ such that 
the following diagram commutes:
\[
\xymatrix{
    A \ar[r]^{f} \ar[d]_{g} & {B} \ar[d]^{g'} \\
    C \ar[r]_{f'} & {D}
  }
\]
and $\forall(y, z)\in B\times C$, if $g'(y)=f'(z)$ then 
$y\in g(A)$ and $z\in h(A)$.\\
We say that $A$ is an 
$[\alpha]$- strong amalgamation basis of $\Gamma$,
if $A$ is a $[\alpha, \alpha]$- strong amalgamation basis. 
\end{df}
In the following remark, We observe that
the most forms of amalgamations given in the 
previous definition, can be characterized by 
the notions of completeness and positive 
completeness which are given in the previous section. 
\begin{rem}
Let $T$ be an h-inductive $L$-theory
and $A$ a model of $T$. We have 
the following properties:
\begin{enumerate}
\item  $A$ is an $h$-amalgamation basis of $T$ if
and only if $T\cup Diag^+(A)$ is positively 
$L(A)$-complete theory.
\item $A$ is an $e$-amalgamation basis of $T$ if
and only if $T\cup Diag(A)$ is positively 
$L(A)$-complete theory.
\item $A$ is an $i$-amalgamation basis of $T$ if
and only if $T\cup Diag^+(A)\cup T_u(A)$ is positively 
$L(A)$-complete theory.
\item 
$A$ is an $[h, e, h, h]$- amalgamation basis of $T$ if
and only if 
the theories $T\cup Diag^+(A)$ and  
$T\cup Diag(A)$ are 
$T$-complete.

\end{enumerate}
\end{rem}
In the following example, we will list some facts on 
amalgamation with the notations and terms given in the 
definition \ref{dfamalgamation}.
\begin{exemple}\label{lemmesamalgamation}
\begin{enumerate}
\item Every $L$-structure $A$ is an $[i,h,s,h]$-amalgamation 
basis in the class of $L$-structures
(lemma 4,
 \cite{ana}). Since every strongly immersion is an 
 immersion, it follows that every 
 $L$-structure $A$ is an $[s,h]$-asymmetric 
 amalgamation basis in the class of $L$-structures.
 \item Every $L$-structure $A$ is an $[s,i]$-asymmetric
  amalgamation basis in the class of $L$-structures
(lemma 5,
 \cite{ana}).
 \item Every $L$-structure $A$ is an $[e,s]$-asymmetric
  amalgamation basis in the class of $L$-structures
(lemma 4,
 \cite{ana3}).
\item Every $L$-structure $A$ is an $[i,h]$-asymmetric
 amalgamation basis in the class of $L$-structures
(lemma 8,
 \cite{begnacpoizat}).
 \item Every pc model of an h-inductive theory $T$ is an 
 $[h]$-amalgamation basis in the class of model of $T$.
 \item The class of $[h]$-amalgamation basis 
 of the theory of rings are the class
of local rings (example 3, \cite{ana2}).
\item Every $L$-structure $A$ is an $[s]$-amalgamation 
basis in the class of $L$-structure. Indeed, suppose
 that $A$ is strongly immersed in two $L$-structures 
 $B$ and $C$. Let $L(B\cup C)$ the language in which 
 we interpret the element of $A$ by the same symbol
 in $B$ and $C$. It is easier to prove that 
 $T_i(C)\cup T_i(B)$ is $L(B\cup C)$-consistent.
\end{enumerate}
\end{exemple}
\begin{thm}\label{strongamalgam1}
Every $L$-structure $A$ is a $[s,i,s,i]$-strong 
amalgamation basis in the class of $L$-structures.
\end{thm}
\preuve Let $A$ be a $L$-structure, 
Suppose that $A$ is immersed 
in a $L$-structure $B$, and strongly immersed in a
$L$-structure $C$.  Suppose that the set 
$$T_i(B)\cup T_u(C)\cup Diag^+(B)\cup Diag^+(C)
\cup \{b\neq c|\ b\in B-A, c\in C-A\}$$
is $L(B\cup C)$-inconsistent. Then there are 
$\neg\psi(\bar a, \bar c)\in T_u(C)$,
$\varphi_1(\bar a, \bar b)\in Diag^+(B)$ and 
$\varphi_2(\bar a, \bar c)\in Diag^+(C)$ such that 
$$T_i(B)\cup \{\neg\psi(\bar a, \bar c),
\varphi_1(\bar a, \bar b), \varphi_2(\bar a, \bar c),
\bigwedge_{i,j}b_i\neq c_j \}$$
is $L(B\cup C)$-inconsistent, thereby  
\begin{equation}\label{eq1}
T_i(B)\vdash\forall\bar y((
\varphi_1(\bar a, \bar b)\wedge
\varphi_2(\bar a, \bar y))\rightarrow 
(\psi(\bar a, \bar y)\vee\bigvee_{i,j}b_i=y_j)).
\end{equation}

Now, since $C\nvDash\psi(\bar a, \bar c)$ and 
$C\vDash\varphi_2(\bar a, \bar c)$, then there is 
$\bar a'\in A$ such that 
$A\nvDash\psi(\bar a, \bar a')$ and 
$A\vDash\varphi_2(\bar a, \bar a')$, 
because otherwise we obtain 
$$A\vdash \forall\bar x (\varphi_2(\bar a, \bar x )
\rightarrow \psi(\bar a, \bar x)).$$
 and given that $C\vdash T_i(A)$, we get a contradiction.\\
So, we obtain  $B\nvDash\psi(\bar a, \bar a')$ and 
$B\vDash\varphi_2(\bar a, \bar a')$. From (\ref{eq1})
we obtain $B\vDash \bigvee_{i,j}b_i=a_j'$.
Since $b_i\not\in A$, we obtain a contradiction.
Then 
$$T_i(B)\cup T_u(C)\cup Diag^+(B)\cup Diag^+(C)
\cup \{b\neq c|\ b\in B-A, c\in C-A\}$$
is $L(B\cup C)$-consistent.\qed 
\begin{thm}\label{strongamalgam2}
Let $T$ be an h-inductive theory.
Every model  $A$ of $T$ is a $[i,i,h, h]$-strong 
amalgamation basis of $T$.
\end{thm}
\preuve 
Let $A$ be a  model of $T$, let $f$ and $g$ 
 two immersion  from $A$ into  $B$
 and $C$ respectively.  We claim that  the set
 $$T\cup Diag^+(B)\cup Diag^+(C)
\cup \{b\neq c|\ b\in B-A, c\in C-A\}$$
is $L(B\cup C)$-consistent. Indeed, otherwise we 
could find $\bar a\in A, \bar b\in B-A, \bar c\in C-A$,
and $\varphi(\bar a, \bar b)\in Diag^+(B),
\psi(\bar a, \bar c)\in Diag^+(C)$ such that 
 $$T\cup\{\varphi(\bar a, \bar b), 
 \psi(\bar a, \bar c), 
 \bigwedge_{i,j}b_i\neq c_j\}$$
 is $L(B\cup C)$-inconsistent. Which implies that 
 $$T\vdash \forall\bar x,\bar y,\bar z\ \ 
 ((\varphi(\bar x, \bar y)\wedge
 \psi(\bar x, \bar y))\rightarrow 
 \bigvee_{i, j}y_i=z_j).$$
 Now, since $C\models\psi(\bar a, \bar c)$ and 
 $A$ is immersed in $C$, then there is 
 $\bar a'\in A$ such that 
 $A\models\psi(\bar a, \bar a')$. Consequently, 
 $B\models\psi(\bar a, \bar a')\wedge 
 \varphi(\bar a, \bar b)$. Given that 
 $T\vdash \forall\bar x,\bar y,\bar z\ \ 
 ((\varphi(\bar x, \bar y)\wedge
 \psi(\bar x, \bar y))\rightarrow 
 \bigvee_{i, j}y_i=z_j)$,
 $\bar b\in B-A$, and $B$ is a model of $B$, 
 we obtain $B\models \bigvee_{i, j}b_i=a'_j$,
 contradiction, which  implies that  
 $T\cup Diag^+(B)\cup Diag^+(C)
\cup \{b\neq c|\ b\in B-A, c\in C-A\}$ 
is consistent. Thereby $A$ is a $[i,i,h, h]$-strong 
amalgamation basis of $T$.\qed
\begin{cor}
Every pc model of an h-inductive theory $T$ is 
an $[h]$-strong amalgamation basis of $T$.
\end{cor}
\begin{lem}\label{lem8}
Every model of an h-inductive theory $T$ is 
an $[i,h,i,h]$-strong amalgamation basis of 
$T$.
\end{lem}
\preuve 
We proceed as in the proof of the theorem \ref{strongamalgam2}. 
Let $A, B$ and $C$ three models
of $T$, such that $A$ is immersed in $B$ and 
continued in $C$ by a homomorphism $f$. The proof
consists of showing the $L(B\cup C)$-consistency
of the set  
$$T_i(C)\cup Diag^+(B)\cup Diag^+(C)
\cup \{b\neq c|\ b\in B-A, c\in C-f(A)\}.$$
Suppose that is not the case, then 
$$T_i(C)\vdash \forall\bar y\ \ 
 ((\varphi(\bar a, \bar c)\wedge
 \psi(\bar a, \bar y))\rightarrow 
 \bigvee_{i, j}y_i=c_j).$$
Given that $B\models\psi(\bar a, \bar b)$ and 
$A$ is immersed $B$, there is $\bar a'\in A$ such that 
$A\models\psi(\bar a, \bar a')$. Which implies
$C\models\varphi(\bar a, \bar c)
\wedge\psi(\bar a, \bar f(\bar a'))$, thereby 
$C\models \bigvee_{i, j}f(\bar a')_i=c_j$,
contradiction.\qed
 \begin{lem}\label{lemma 9}
 Let $B$ be a $[h]$-strong  amalgamation basis
  of $T$ and  $A$ a model of $T$ that is 
  immersed in $B$, then  $A$ is a $[h]$-strong
   amalgamation basis of $T$. 
 \end{lem}
\preuve 
 Suppose that $A$ is continued into $C$ and 
 $D$ two models of $T$ by $f_1$ and 
 $f_2$ respectively. Given that $A$ is immersed 
 in $B$, and every $L$-structure  is an 
 $[i,h,s,h]$-strong amalgamation basis in the class of models of $T$
 (lemma \ref{lem8}), we obtain the 
 following commutative diagram:
 $$
  \xymatrix{
  &C \ar[r]^{i_1}&C'\\
    A \ar[r]^{i} \ar[rd]_{g_1} \ar[ru]^{f_1}&
   {B} \ar[rd]_{g_2} \ar[ru]^{f_2} \\
    & D \ar[r]_{i_2} & {D'}
  }
  $$
 where $i_1$ and $i_2$ are immersions, 
 $f_2$ and $g_2$ homomorphisms and 
 $C', D'$ two models of $T$. Now since $B$ is 
 a $[h]$-strong  amalgamation basis
  of $T$, we complete the previous  diagram and 
  we get the following commutative diagram:
  $$
  \xymatrix{
  &C \ar[r]^{i_1}&C'\ar[rd]^{f_3}&&\\
    A \ar[r]^{i} \ar[rd]_{g_1} \ar[ru]^{f_1}&
   {B} \ar[rd]_{g_2} \ar[ru]^{f_2}&& E \\
    & D \ar[r]_{i_2} & {D'}\ar[ru]_{g_3}&&
  }
  $$
  where $E$ is a model of $T$, $f_3$ and $g_3$ homomorphisms such that 
  $$\forall c\in C'-f_2(B), \forall d\in D'-g_2(B);\ \ \  f_3(c)\neq 
 g_3(d).$$
 Given that the diagram is commutative.
  $A$, $C$ and $D$ are immersed in 
 $B$, $C'$ and $D'$ repetitively. Then 
 $$\forall c\in C-f_1(A), \forall d\in D-g_2(A);\ \ \  s_1\circ f_3(c)\neq 
 s_2\circ g_3(d).$$

 \bibliographystyle{plain}

\end{document}